\input amstex
\input amsppt.sty
\magnification=\magstep1 \hsize=36truecc
 \vsize=23.5truecm
\baselineskip=14truept
 \NoBlackBoxes
\def\q{\quad}
\def\qq{\qquad}
\def\mod#1{\ (\text{\rm mod}\ #1)}

\def\t{\text}
\def\qtq#1{\q\t{#1}\q}
\par Preprint: July 6, 2010
\def\mod#1{\ (\text{\rm mod}\ #1)}
\def\qtq#1{\q\t{#1}\q}
\def\f{\frac}
\def\e{\equiv}
\def\b{\binom}
\def\la{\lambda}

 \def\ls#1#2{\big(\f{#1}{#2}\big)}

\let \pro=\proclaim
\let \endpro=\endproclaim
\topmatter
\title Constructing $x^2$ for primes $p=ax^2+by^2$
\endtitle
\author ZHI-Hong Sun\endauthor
\affil School of the Mathematical Sciences, Huaiyin Normal
University,
\\ Huaian, Jiangsu 223001, PR China
\\ E-mail: szh6174$\@$yahoo.com
\\ Homepage: http://www.hytc.edu.cn/xsjl/szh
\endaffil

 \nologo \NoRunningHeads

\abstract{Let $a$ and $b$ be positive integers and let $p$ be an odd
prime such that $p=ax^2+by^2$ for some integers $x$ and $y$. Let
$\lambda(a,b;n)$ be given by
$q\prod_{k=1}^{\infty}(1-q^{ak})^3(1-q^{bk})^3=\sum_{n=1}^{\infty}\lambda(a,b;n)q^n.$
In the paper, using Jacobi's identity
$\prod_{n=1}^{\infty}(1-q^n)^3=\sum_{k=0}^{\infty}(-1)^k(2k+1)q^{\frac{k(k+1)}2}$
we construct $x^2$ in terms of $\lambda(a,b;n)$. For example, if
$2\nmid ab$ and $p\nmid ab(ab+1)$, then
$(-1)^{\frac{a+b}2x+\frac{b+1}2}(4ax^2-2p)=\lambda(a,b;((ab+1)p-a-b)/8+1)$.
We also give formulas for $\lambda(1,3;n+1),\lambda(1,7;2n+1)$,
$\lambda(3,5;2n+1)$ and $\lambda(1,15;4n+1)$.
\par\q
\newline MSC: Primary 11E16, Secondary 11E25 \newline Keywords:
Binary quadratic form; Jacobi's identity}
 \endabstract
  \footnote"" {The author is
supported by the Natural Sciences Foundation of China (grant No.
10971078).}
\endtopmatter
\document

\subheading{1. Introduction}\par\q

\par Let $p$ be a prime of the form $4k+1$. The two squares theorem asserts
 that there are unique
 positive integers $x$ and $y$ such that $p=x^2+y^2$ and $2\nmid x$. Since
Legendre and Gauss, there are several methods to construct $x$ and
$y$. For example, if we choose the sign of $x$ so that $x\e 1\mod
4$, we then have
$$\align &(1.1)\q \t{(Gauss[3],1825)}\qq 2x\e \b{\f{p-1}2}{\f{p-1}4}\mod p,
\\&(1.2)\q \t{(Jacobsthal[3],1907)}\qq 2x=-\sum_{n=0}^{p-1}\ls {n^3-4n}p,
\\&(1.3)\q  \t{(Liouville[7], 1862)}\qq 6x=N(p=t^2+u^2+v^2+16w^2)-3p-3,
\\&(1.4)\q  \t{(Klein and Fricke[6],1892)}\qq
4x^2-2p=[q^p]q\prod_{k=1}^{\infty}(1-q^{4k})^6,
\\&(1.5)\q \t{(Sun[13], 2006)}\qq 2y=5p+3-8V_p(z^4-3z^2+2z)\qtq{for}p\e 5\mod {12},\endalign$$
where $\ls ap$ is the Legendre-Jacobi-Kronecker symbol,
$N(p=t^2+u^2+v^2+16w^2)$ is the number of integral solutions to
$p=t^2+u^2+v^2+16w^2$, $[q^n]f(q)$ denotes the coefficient of $q^n$
in the power series expansion of $f(q)$, and $V_p(f(z))$ is the
number of $c\in\{0,1,\ldots,p-1\}$ such that $f(z)\e c\mod p$ is
solvable. We note that (1.3) was conjectured by Liouville  and
proved by A. Alaca, S. Alaca, M. F. Lemire, and K. S. Williams
([1]).

\par Let $\Bbb Z$ and $\Bbb N$ be the sets of integers and positive
integers, respectively. For $a,b,n\in\Bbb N$ let $\lambda(a,b;n)$ be
given by
$$q\prod_{k=1}^{\infty}(1-q^{ak})^3(1-q^{bk})^3=\sum_{n=1}^{\infty}\lambda(a,b;n)q^n
\q (|q|<1).$$ In his lost notebook, Ramanujan ([9]) conjectured that
$\lambda(1,7;n)$ is multiplicative and
$$\align&\sum_{n=1}^{\infty}\f{\lambda(1,7;n)}{n^s}\\&=\f 1{1+7^{1-s}}\prod_{p\e
3,5,6\mod 7}\f 1{1-p^{2-2s}}\prod_{p\e 1,2,4\mod 7}\f
1{1-(4x^2-2p)p^{-s}+p^{2-2s}},\endalign$$ where $s>1$, $p$ runs over
all distinct primes and $x^2$ is given by  $p=x^2+7y^2\e 1,2,4\mod
7$. This was proved by Hecke[5]. See also [10]. The above assertion
of Ramanujan implies
$$\lambda(1,7;p)=4x^2-2p\qtq{for primes}p=x^2+7y^2\e 1,2,4\mod 7.\tag 1.6$$
In his lost notebook, Ramanujan[9] also conjectured that
$\la(4,4;n)$ is multiplicative. This was proved by Mordell[8] in
1917. It is easily seen that $\la(4,4;p)=\la(1,1;(p+3)/4)$ for $p\e
1\mod 4$. Thus, (1.4) is equivalent to
  $$\la(1,1;(p+3)/4)=4x^2-2p\qtq{for primes} p=x^2+y^2\e 1\mod
4\ \t{with}\ 2\nmid x.\tag 1.7$$   In 1985 Stienstra and Beukers[11]
proved
$$\la(2,6;p)=4x^2-2p\qtq{for primes}p=x^2+3y^2\e 1\mod 3.\tag 1.8$$
It is easily seen that $\la(2,6;p)=\la(1,3;(p+1)/2)$ for odd $p$.
\par
 In the paper, with the help of the following Jacobi's identity
 ([2])
$$\prod_{n=1}^{\infty}(1-q^n)^3=\sum_{k=0}^{\infty}(-1)^k(2k+1)q^{\f{k(k+1)}2}
\q(|q|<1),\tag 1.9$$ we construct $x^2$ for primes $p=ax^2+by^2$.
For example, if  $a,b\in\Bbb N$, $2\nmid ab$ and $p$ is an odd prime
such that $p\nmid ab(ab+1)$ and $p=ax^2+by^2$ with $x,y\in\Bbb Z$,
then
$$\aligned&(-1)^{\f{a+b}2x+\f{b+1}2}(4ax^2-2p)
=\lambda(a,b;n+1)\\&=\sum_{k_1+2k_2+\cdots+nk_n=n}
 (-3)^{k_1+\cdots+k_n}\f{(a\sigma(\f 1a)+b\sigma(\f 1b))^{k_1}
 \cdots(a\sigma(\f na)+b\sigma(\f nb))^{k_n}}
 {1^{k_1}\cdot k_1!\cdots n^{k_n}\cdot k_n!},\endaligned\tag 1.10$$
where $n=((ab+1)p-a-b)/8$ and $$\sigma(m)=\cases\sum\limits_{d\mid
m}d&\t{if $m\in\Bbb N$,}
\\0&\t{otherwise.}\endcases\tag 1.11$$ This can be viewed as a vast
generalization of (1.6)-(1.8). In the paper we also give formulas
for $\la(1,3;n+1),\la(1,7;2n+1)$, $\la(3,5;2n+1)$ and
$\la(1,15;4n+1)$.
 \subheading{2. Basic lemmas}
\par A negative integer $d$ with $d\e 0,1\mod 4$ is called a discriminant.
 Let $d$ be a discriminant. The conductor of $d$ is the largest positive
 integer $f=f(d)$ such that $d/f^2\e 0,1\mod 4$. As usual we set $w(d)=2,4,6$
according as $d<-4,d=-4$ or $d=-3$. For $a,b,c\in\Bbb Z$ we denote
the equivalence class containing the form $ax^2+bxy+cy^2$ by
$[a,b,c]$.
 Let $H(d)$ be
the form class group consisting of classes of primitive, integral
binary quadratic forms of discriminant $d$. See [4]. For $n\in\Bbb
N$ and $[a,b,c]\in H(d)$, following [14] we define
$$R([a,b,c],n)=|\{\langle x,y\rangle\in\Bbb Z\times\Bbb Z:\ n=ax^2+bxy+cy^2\}|.$$
It is known that $R([a,b,c],n)=R([a,-b,c],n)$. If $R([a,b,c],n)>0$,
we say that $n$ is represented by $[a,b,c]$.
\par For $m,n\in\Bbb N$ let $(m,n)$ be the greatest common divisor
of $m$ and $n$. \pro{Lemma 2.1 ([14, Lemma 5.2])} Let $d<0$ be a
discriminant with conductor $f$. Let $p$ be a prime and $K\in H(d)$.
\par $(\t{\rm i})$ $p$ is
represented by some class in $H(d)$ if and only if $(\f dp)=0,1$ and
$p\nmid f$.
\par $(\t{\rm ii})$ Suppose $p\mid d$ and $p\nmid f$. Then $p$ is
represented by exactly one class $A\in H(d)$, and $A=A^{-1}$.
Moreover, $R(A,p)=w(d)$.
\par $(\t{\rm iii})$ Suppose $(\f dp)=1$. Then $p$ is represented
by some class $A\in H(d)$, and
$$R(K,p)=\cases 0&\t{if $K\not=A,A^{-1}$,}
\\ w(d)&\t{if $A\not=A^{-1}$ and $K\in\{A,A^{-1}\}$,}
\\ 2w(d)&\t{if $K=A=A^{-1}$.}\endcases$$\endpro

\pro{Lemma 2.2 ([14, Theorem 7.1])} Let $d$ be a negative
discriminant and $K\in H(d)$. If $n_1,n_2\in\Bbb N$ and
$(n_1,n_2)=1$, then
$$R(K,n_1n_2)=\f 1{w(d)}\sum\Sb K_1K_2=K\\K_1,K_2\in H(d)\endSb
R(K_1,n_1)R(K_2,n_2).$$
\endpro

\pro{Lemma 2.3} Let $a,b\in\Bbb N$ and let $p$ be an odd prime such
that $p\not=a,b$, $p\nmid ab+1$ and $p=ax^2+by^2$ with  $x,y\in\Bbb
Z$. \par $(\t{\rm i})$ If $ab+1$ is not a square, then
$R([a,0,b],(ab+1)p)=8$ and all the integral solutions to the
equation $(ab+1)p=aX^2+bY^2$ are given by $\{x\pm by,ax\mp
y\},\{x\pm by,-(ax\mp y)\},\{-(x\pm by),ax\mp y\}$ and $\{-(x\pm
by),-(ax\mp y)\}.$
\par $(\t{\rm ii})$ If $ab+1=m^2$  for $m\in\Bbb N$, then
$R([a,0,b],(ab+1)p)=12$ and all the integral solutions to the
equation $(ab+1)p=aX^2+bY^2$ are given by $\{mx,\pm my\},\{-mx,\pm
my\}$, $\{x\pm by,ax\mp y\},\{x\pm by,-(ax\mp y)\},\{-(x\pm
by),ax\mp y\}$ and $\{-(x\pm by),-(ax\mp y)\}.$

\endpro
Proof. Since $p\not=a,b$ and $p=ax^2+by^2$, we see that $p\nmid ab$
and $(a,b)=1$. As $p\nmid ab+1$ and $[1,0,ab][a,0,b]=[a,0,b]$, by
Lemmas 2.1 and 2.2 we have
$$\aligned R([a,0,b],(ab+1)p)&=\f 1{w(-4ab)}\sum\Sb AB=[a,0,b]\\A,B\in
H(-4ab)\endSb R(A,p)R(B, ab+1)
\\&=\f{R([a,0,b],p)R([1,0,ab],ab+1)}{w(-4ab)}
=2R([1,0,ab],ab+1).\endaligned$$ If $ab+1$ is not a square and
$ab+1=X^2+abY^2$ for some $X,Y\in\Bbb Z$, we must have $X^2=Y^2=1$
 and so $R([1,0,ab],ab+1)=4$. Hence
$R([a,0,b],(ab+1)p)=2R([1,0,ab],ab+1)=8$. It is clear that
$$xy\not=0\qtq{and} (ab+1)p=(ab+1)(ax^2+by^2)=a(x\pm by)^2+b(ax\mp y)^2.$$ Thus,
$\{x\pm by,ax\mp y\},\{x\pm by,-(ax\mp y)\},\{-(x\pm by),ax\mp
y\},\{-(x\pm by),-(ax\mp y)\}$ are the eight integral solutions to
the equation $(ab+1)p=aX^2+bY^2$. This proves (i).
\par If $ab+1=m^2$ for $m\in\Bbb N$ and
$ab+1=X^2+abY^2$ for some $X,Y\in\Bbb Z$, we must have $Y\in\{0,\pm
1\}$  and so $R([1,0,ab],ab+1)=6$. Hence
$R([a,0,b],(ab+1)p)=2R([1,0,ab],ab+1)=12$. Since $xy\not=0$ and
$$(ab+1)p=(ab+1)(ax^2+by^2)=a(mx)^2+b(my)^2=a(x\pm by)^2+b(ax\mp
y)^2,$$ we see that $\{mx,\pm my\},\{-mx,\pm my\},\{x\pm by,ax\mp
y\},\{x\pm by,-(ax\mp y)\},\{-(x\pm by),ax\mp y\},\{-(x\pm
by),-(ax\mp y)\}$ are $12$ integral solutions to the equation
$(ab+1)p=aX^2+bY^2$. This proves (ii).

\pro{Lemma 2.4} Let $a,b\in\Bbb N$, $(a,b)=1$ and let $p$ be an odd
prime such that $p\not=ab,ab+1$ and $p=x^2+aby^2$ with $x,y\in\Bbb
Z$. Suppose $(a-1)(b-1)\not=0$ or $a+b$ is not a square. Then
$R([a,0,b],(a+b)p)=8$ and all the integral solutions to the equation
$(a+b)p=aX^2+bY^2$ are given by
$$\{x\pm by,x\mp ay\},\{x\pm by,-(x\mp ay)\},\{-(x\pm by),x\mp
ay\},\{-(x\pm by),-(x\mp ay)\}.$$
\endpro
Proof. Since $p\not=ab,ab+1$, we see that $p=x^2+aby^2>1+ab\ge a+b$
and so $p\nmid a+b$.
 As  $[1,0,ab][a,0,b]=[a,0,b]$, by Lemmas 2.1 and
2.2 we have
$$\aligned R([a,0,b],(a+b)p)&=\f 1{w(-4ab)}\sum\Sb AB=[a,0,b]\\A,B\in
H(-4ab)\endSb R(A,p)R(B, a+b)
\\&=\f 1{w(-4ab)}R([1,0,ab],p)R([a,0,b],a+b)
=2R([a,0,b],a+b).\endaligned$$ If $a+b=aX^2+bY^2$ for some
$X,Y\in\Bbb Z$, we must have $X^2=Y^2=1$. Thus $R([a,0,b],a+b)=4$
and so $R([a,0,b],(a+b)p)=2R([a,0,b],a+b)=8$. It is clear that
$$xy\not=0\qtq{and} (a+b)p=(a+b)(x^2+aby^2)=a(x\pm by)^2+b(x\mp ay)^2.$$ Thus,
$\{x\pm by,x\mp ay\},\{x\pm by,-(x\mp ay)\},\{-(x\pm by),x\mp
ay\},\{-(x\pm by),-(x\mp ay)\}$ are the eight integral solutions to
the equation $(a+b)p=aX^2+bY^2$. This completes the proof.
 \pro{Lemma 2.5} Let
$a,b,n\in\Bbb N$. Then
$$\sum\Sb x,y\in\Bbb Z,x\e y\e 1\mod
4\\ax^2+by^2=8n+a+b\endSb xy=\lambda(a,b;n+1).$$
\endpro
Proof. Using Jacobi's identity (1.9) we see that
$$\aligned&q\prod_{n=1}^{\infty}(1-q^{an})^3(1-q^{bn})^3
\\&=q\Big(\sum_{k=0}^{\infty}(-1)^k(2k+1)q^{a\f{k(k+1)}2}\Big)
\Big(\sum_{m=0}^{\infty}(-1)^m(2m+1)q^{b\f{m(m+1)}2}\Big)
\\&=\sum_{n=0}^{\infty}\sum\Sb k,m\ge 0\\a\f{k(k+1)}2+b\f{m(m+1)}2=n
\endSb (-1)^k(2k+1)\cdot (-1)^m(2m+1)q^{n+1}.\endaligned$$
Thus,
$$\aligned \la(a,b;n+1)&=\sum\Sb k,m\ge 0\\a\f{k(k+1)}2+b\f{m(m+1)}2=n
\endSb (-1)^k(2k+1)\cdot (-1)^m(2m+1)
\\&=\sum\Sb k,m\ge 0\\a(2k+1)^2+b(2m+1)^2=8n+a+b
\endSb (-1)^k(2k+1)\cdot (-1)^m(2m+1)
\\&=\sum\Sb x\e y\e 1\mod 4\\ax^2+by^2=8n+a+b\endSb xy.\endaligned$$
This proves the lemma.

\pro{Lemma 2.6} Let $a,b\in\Bbb N$ with $(a,b)=1$ and $ab\e 1\mod
4$. Let $p$ be an odd prime such that $R([a,0,b],2p)>0$. Then
$R([a,0,b],2p)=2w(-4ab)$.
\endpro
Proof. Suppose that $2p=ax^2+by^2$ with $x,y\in\Bbb Z$. We claim
that $p\nmid ab$. If $p\mid a$, then  $p\mid by^2$ and so $p\mid b$.
This contradicts the fact $(a,b)=1$. Hence $p\nmid a$. Similarly, we
have $p\nmid b$. Since $-ab\e 3\mod 4$ we see that $2\nmid f(-4ab)$.
Thus, by Lemma 2.1, there exists exactly one class $A\in H(-4ab)$
such that $R(A,2)>0$ and we have $A=A^{-1}$. Using Lemmas 2.1, 2.2
and the fact $R([a,0,b],2p)>0$ we see that $R([a,0,b],2p)=\f
1{w(-4ab)}R(A,2)R(A[a,0,b],p)=R(A[a,0,b],p)=2w(-4ab)$. This
completes the proof.
\par\q

 \pro{Lemma 2.7} Let $a,b\in\Bbb N$, $ab\e
3\mod 4$, $K\in H(-4ab)$ and $K=K^{-1}$. Let $p$ be an odd prime
such that $p\nmid ab$ and $R(K,4p)>0$. Then $R(K,4p)=2w(-ab)$.
\endpro
Proof. From Lemma 2.2 we have
$$2R(K,4p)=\sum\Sb AB=K\\A,B\in H(-4ab)\endSb R(A,p)R(B,4)>0.$$Since
$2\mid f(-4ab)$, by [14, Theorem 5.3(i)] we have $R(B,4)=0$ or
$w(-ab)$ for $B\in H(-4ab)$. Suppose $R(A,p)>0$ for $A\in H(-4ab)$.
 Then $(AK)^{-1}=K^{-1}A^{-1}=KA^{-1}=A^{-1}K$ and so
 $R(A^{-1}K,4)=R((AK)^{-1},4)=R(AK,4)$.
From the above and Lemma 2.1 we see that
$$\aligned&2R(K,4p)\\&=\cases R(A,p)R(AK,4)=4\cdot w(-ab)&\t{if $A=A^{-1}$,}
\\R(A,p)R(A^{-1}K,4)+R(A^{-1},p)R(AK,4)=2 w(-ab)+2 w(-ab)
&\t{if $A\not=A^{-1}$.}
\endcases\endaligned$$ This yields the result.
\par\q
\par  If $\{a_n\}$ and $\{b_n\}$ are two
sequences satisfying
$$a_1=b_1\qtq{and}b_n+a_1b_{n-1}+\cdots+a_{n-1}b_1=na_n\ (n=2,3,\ldots),$$
we say that $(a_n,b_n)$ is a Newton-Euler pair as in [12].  For a
rational number $m$ let $\sigma(m)$ be given by (1.11). Now we state
the following result.
 \pro{Lemma 2.8 } Let $a,b\in\Bbb N$. Then
$(\lambda(a,b;n+1),-3(a\sigma(n/a)+b\sigma(n/b))$ is a Newton-Euler
pair. That is, for $n\in\Bbb N$,
$$a\sigma(\f na)+b\sigma(\f nb)+\sum_{k=1}^{n-1}
\big(a\sigma(\f ka)+b\sigma(\f kb)\big)\lambda(a,b;n+1-k) =-\f
n3\lambda(a,b;n+1).$$
\endpro
Proof. Suppose $q\in\Bbb R$ and $|q|<1$. As
$$1-q^n=\prod_{r=0}^{n-1}\big(1-\t{e}^{2\pi i\f rn}q\big),$$
we see that
$$\align 1+\sum_{n=1}^{\infty}\la(a,b;n+1)q^n
&=\prod_{k=1}^{\infty}\big(1-q^{ak}\big)^3\big(1-q^{bk}\big)^3
\\&=\prod_{k=1}^{\infty}\prod_{r=0}^{ak-1}\big(1-\t{e}^{2\pi i\f r{ak}}q
\big)^3 \prod_{s=0}^{bk-1}\big(1-\t{e}^{2\pi i\f s{bk}}q\big)^3.
 \endalign$$
 Observe that
 $$\aligned &\sum_{k=1}^{\infty}\Big\{\sum_{r=0}^{ak-1}3\Big(\t{e}^{2\pi i\f
 r{ak}}\Big)^n+\sum_{s=0}^{bk-1}3\Big(\t{e}^{2\pi i\f
 s{bk}}\Big)^n\Big\}
 \\&=3\sum\Sb k\in\Bbb N\\ak\mid n\endSb ak
 +3\sum\Sb k\in\Bbb N\\bk\mid n\endSb bk=3a\sigma\big(\f na\big)+3b\sigma\big(\f nb\big).
 \endaligned$$
 From the above and [12, Example 1, p.103] we deduce the result.
 \pro{Lemma 2.9} Let $a,b,n\in\Bbb
N$. Then
 $$\align&\la(a,b;n+1)\\&=\sum_{k_1+2k_2+\cdots+nk_n=n}
 (-3)^{k_1+\cdots+k_n}\f{(a\sigma(\f 1a)+b\sigma(\f 1b))^{k_1}
 \cdots(a\sigma(\f na)+b\sigma(\f nb))^{k_n}}
 {1^{k_1}\cdot k_1!\cdots n^{k_n}\cdot k_n!}.\endalign$$
 \endpro
 Proof. This is immediate from Lemma 2.8 and [12, Theorem 2.2].

\subheading{3. Constructing $x^2$ for primes $p=ax^2+by^2$}
 \par\q
\pro{Theorem 3.1} Let $a,b\in\Bbb N$ with $2\nmid ab$.
 Let $p$ be an odd prime such that
$p\not=a,b$, $p\nmid ab+1$ and $p=ax^2+by^2$ with $x,y\in\Bbb Z$.
Let $n=((ab+1)p-a-b)/8$. Then
$$\align &(-1)^{\f{a+b}2x+\f{b+1}2}(4ax^2-2p)
\\&=\lambda(a,b;n+1)=\sum_{k_1+2k_2+\cdots+nk_n=n}
 (-3)^{k_1+\cdots+k_n}\f{(a\sigma(\f 1a)+b\sigma(\f 1b))^{k_1}
 \cdots(a\sigma(\f na)+b\sigma(\f nb))^{k_n}}
 {1^{k_1}\cdot k_1!\cdots n^{k_n}\cdot k_n!}.\endalign$$
\endpro
Proof. Clearly $2\mid x$ or $2\mid y$. If $2\mid y$, then $p\e
ax^2\e a\mod 4$ and so $(ab+1)p\e (ab+1)a\e a+b\mod 8$. If $2\mid
x$, then $p\e by^2\e b\mod 4$ and so $(ab+1)p\e (ab+1)b\e a+b\mod
8$. Thus $n\in\Bbb N$. By Lemma 2.3, all the integral solutions
$\{X,Y\}$ with $2\nmid XY$ to the equation
$8n+a+b=(ab+1)p=aX^2+bY^2$ are given by $\{x\pm by,ax\mp y\},\{x\pm
by,-(ax\mp y)\},\{-(x\pm by),ax\mp y\},\{-(x\pm by),-(ax\mp y)\}$.
 Since  $x\pm by\e (-1)^{\f{a+b}2x+\f{b+1}2}(ax\mp y)
 \mod 4$, applying Lemma 2.5 we have
$$\aligned\la (a,b;n+1)&=\sum\Sb X\e Y\e 1\mod 4\\aX^2+bY^2=8n+a+b\endSb XY
\\&=(x+by)\cdot (-1)^{\f{a+b}2x+\f{b+1}2}(ax-y)+(x-by)\cdot
(-1)^{\f{a+b}2x+\f{b+1}2}(ax+y)\\&=(-1)^{\f{a+b}2x+\f{b+1}2}2(ax^2-by^2)
=(-1)^{\f{a+b}2x+\f{b+1}2}(4ax^2-2p).\endaligned$$ This together
with Lemma 2.9 yields the result.
 \pro{Corollary 3.1} Let $p$ be a prime of
the form $4k+1$ and so $p=x^2+y^2$ with $x,y\in\Bbb Z$ and $2\nmid
x$. Let $n=(p-1)/4$. Then $$4x^2-2p=\sum_{k_1+2k_2+\cdots+nk_n=n}
 (-6)^{k_1+\cdots+k_n}\f{\sigma(1)^{k_1}
 \cdots\sigma(n)^{k_n}}
 {1^{k_1}\cdot k_1!\cdots n^{k_n}\cdot k_n!}.$$
\endpro
Proof. Taking $a=b=1$ in Theorem 3.1 we obtain the result.

\pro{Corollary 3.2} Suppose that $p\e 1,9\mod{20}$ is a prime and so
$p=x^2+5y^2$ for some $x,y\in\Bbb Z$. Let $n=3(p-1)/4$. Then
$$\align &(-1)^{x-1}(4x^2-2p)=\la(1,5;(3p+1)/4)
\\&=\sum_{k_1+2k_2+\cdots+nk_n=n}
 (-3)^{k_1+\cdots+k_n}\f{(\sigma(1)+5\sigma(\f 15))^{k_1}
 \cdots(\sigma(n)+5\sigma(\f n5))^{k_n}}
 {1^{k_1}\cdot k_1!\cdots n^{k_n}\cdot k_n!}.\endalign$$
\endpro
Proof. Taking $a=1$ and $b=5$ in Theorem 3.1 we obtain the result.
 \pro{Theorem 3.2} Let $a,b\in\Bbb N$ with $(a,b)=1$. Let $p$ be an odd prime such that
$p\not=ab,ab+1$ and $p=x^2+aby^2$ with $x,y\in\Bbb Z$. Let
$n=(a+b)(p-1)/8$.
\par $(\t{\rm i})$ If $2\nmid ab$, then
$$\align&(-1)^{\f{ab+1}2y}(4x^2-2p)=\lambda(a,b;n+1)\\&=\sum_{k_1+2k_2+\cdots+nk_n=n}
 (-3)^{k_1+\cdots+k_n}\f{(a\sigma(\f 1a)+b\sigma(\f 1b))^{k_1}
 \cdots(a\sigma(\f na)+b\sigma(\f nb))^{k_n}}
 {1^{k_1}\cdot k_1!\cdots n^{k_n}\cdot k_n!}.\endalign$$
\par $(\t{\rm ii})$ If $2\nmid a$, $2\mid b$, $8\nmid b$ and $8\mid p-1$, then
$$\align&(-1)^{\f y2}(4x^2-2p)=\lambda(a,b;n+1)\\&=\sum_{k_1+2k_2+\cdots+nk_n=n}
 (-3)^{k_1+\cdots+k_n}\f{(a\sigma(\f 1a)+b\sigma(\f 1b))^{k_1}
 \cdots(a\sigma(\f na)+b\sigma(\f nb))^{k_n}}
 {1^{k_1}\cdot k_1!\cdots n^{k_n}\cdot k_n!}.\endalign$$
\endpro
Proof. If $(a-1)(b-1)\not=0$ or $a+b$ is not a square, using Lemma
2.4 we see that $R([a,0,b],(a+b)p)=8$ and all the integral solutions
to $(a+b)p=aX^2+bY^2$ are given by $\{x\pm by,x\mp ay\},\{x\pm
by,-(x\mp ay)\},\{-(x\pm by),x\mp ay\},\{-(x\pm by),-(x\mp ay)\}.$
If $a=1$ and $b+1=m^2$ for $m\in\Bbb N$, using Lemma 2.3(ii) we see
that $R([1,0,b],(b+1)p)=12$ and all the integral solutions to
$(b+1)p=X^2+bY^2$ are given by $\{mx,\pm my\},\{-mx,\pm my\},\{x\pm
by,x\mp y\},\{x\pm by,-(x\mp y)\},\{-(x\pm by),x\mp y\},\{-(x\pm
by),-(x\mp y)\}.$ If $b=1$ and $a+1=k^2$ for $k\in\Bbb N$, using
Lemma 2.3(ii) we see that $R([a,0,1],(a+1)p)=12$ and all the
integral solutions to $(a+1)p=aX^2+Y^2$ are given by $\{ky,\pm
kx\},\{-ky,\pm kx\},\{x\pm y,x\mp ay\},\{x\pm y,-(x\mp
ay)\},\{-(x\pm y),x\mp ay\},\{-(x\pm y),-(x\mp ay)\}.$
\par We first
assume $2\nmid ab$. If $ab\e 1\mod 4$, then $p=x^2+aby^2\e 1\mod 4$
and so $(a+b)(p-1)\e 0\mod 8$. If $ab\e 3\mod 4$, then $4\mid a+b$
and so $8\mid (a+b)(p-1)$. Thus, we always have $8\mid (a+b)(p-1)$.
  It is easily seen that $x\pm by\e 1\mod 2$ and $x\pm by\e (-1)^{\f{ab+1}2y}(x\mp ay)
 \mod 4$. Thus, applying the above and Lemma 2.5 we have
$$\aligned\la (a,b;n+1)&=\sum\Sb X\e Y\e 1\mod 4\\aX^2+bY^2=8n+a+b\endSb XY
\\&=(x+by)\cdot (-1)^{\f{ab+1}2y}(x-ay)+(x-by)\cdot
(-1)^{\f{ab+1}2y}(x+ay)\\&=(-1)^{\f{ab+1}2y}2(x^2-aby^2)=(-1)^{\f{ab+1}2y}(4x^2-2p).
\endaligned$$  This together with Lemma 2.9 proves (i).

\par Now we consider (ii). Since $2\nmid a$, $2\mid b$, $8\nmid b$ and $8\mid p-1$, we deduce
$2\nmid x$, $8\mid by^2$ and so $2\mid y$.  It is easily seen that
$x\pm by\e 1\mod 2$ and $x\pm by\e (-1)^{\f y2}(x\mp ay)
 \mod 4$. Thus, applying the above and Lemma 2.5 we have
$$\aligned\la (a,b;n+1)&=\sum\Sb X\e Y\e 1\mod 4\\aX^2+bY^2=8n+a+b\endSb XY
\\&=(x+by)\cdot (-1)^{\f y2}(x-ay)+(x-by)\cdot
(-1)^{\f y2}(x+ay)\\&=(-1)^{\f y2}2(x^2-aby^2)=(-1)^{\f
y2}(4x^2-2p).
\endaligned$$  This together with Lemma 2.9 yields (ii). The proof is now complete.

\pro{Corollary 3.3} Let $a,b\in\Bbb N$ with $2\nmid ab$ and
$(a,b)=1$. Let $p$ be an odd prime such that $p\not=ab,ab+1$ and
$p=x^2+aby^2$ with $x,y\in\Bbb Z$. Then
$$\lambda\Big(a,b;\f{(a+b)(p-1)}8+1\Big)=\lambda\Big(1,ab;\f{(ab+1)(p-1)}8+1\Big).$$
\endpro
Proof. By Theorem 3.1 we have
$$(-1)^{\f{1+ab}2(x+1)}(4x^2-2p)=\lambda\Big(1,ab;\f{(ab+1)(p-1)}8+1\Big).$$
This together with Theorem 3.2(i) gives the result.

\pro{Corollary 3.4} Suppose $a\in\Bbb N$ and $2\nmid a$. Let $p$ be
an odd  prime such that  $p=x^2+16ay^2$ with $x,y\in\Bbb Z$. Then
$$(-1)^y(4x^2-2p)
=\lambda\Big(a,4;\f{(a+4)p-a+4}8\Big).$$
\endpro
Proof. Taking $b=4$ and replacing $y$ with $2y$ in Theorem 3.2(ii)
we deduce the result.
\par\q
\par Let $p$ be an odd prime. From Theorem 3.2 we deduce:
\
$$\align &(-1)^{\f
y2}(4x^2-2p)=\la(1,2;(3p+5)/8)\qtq{for}p=x^2+2y^2\e 1\mod 8,\tag 3.1
\\&(-1)^{\f
y2}(4x^2-2p)=\la(1,6;(7p+1)/8)\qtq{for}p=x^2+6y^2\e 1\mod {24}, \tag
3.2\\&(-1)^{\f y2}(4x^2-2p)=\la(1,10;(11p-3)/8)\q \t{for}\q
p=x^2+10y^2\e 1,9\mod {40},\tag 3.3
\\&(-1)^{\f
y2}(4x^2-2p)=\la(1,12;(13p-5)/8)\q \t{for}\q p=x^2+12y^2\e 1\mod
{24}.\tag 3.4
\endalign$$

 \pro{Theorem 3.3} Let $a,b\in\Bbb N$,
$2\nmid a$, $2\mid b$ and $8\nmid b$. Let $p$ be a prime such that
$p\e a\mod 8$, $p\not=a$, $p\nmid ab+1$ and $p=ax^2+by^2$ with
$x,y\in\Bbb Z$. Let $n=((ab+1)p-a-b)/8$. Then
$$\align&(-1)^{\f{a-1}2+\f y2}(4ax^2-2p)=(-1)^{\f{a-1}2+\f y2}(2p-4by^2)
=\lambda(a,b;n+1)
\\&=\sum_{k_1+2k_2+\cdots+nk_n=n}
 (-3)^{k_1+\cdots+k_n}\f{(a\sigma(\f 1a)+b\sigma(\f 1b))^{k_1}
 \cdots(a\sigma(\f na)+b\sigma(\f nb))^{k_n}}
 {1^{k_1}\cdot k_1!\cdots n^{k_n}\cdot k_n!}.\endalign$$
\endpro
Proof. Clearly we have $2\nmid x$ and so $8\mid by^2$. Since $8\nmid
b$ we must have $2\mid y$. As $p\e a\mod 8$ we have $(ab+1)p\e
(1+ab)a\e a+b\mod 8$. Set $n=\f{(ab+1)p-a-b}8$. By Lemma 2.3, all
the integral solutions $\{X,Y\}$ with $2\nmid XY$ to the equation
$8n+a+b=(ab+1)p=aX^2+bY^2$ are given by $\{x\pm by,ax\mp y\},\{x\pm
by,-(ax\mp y)\},\{-(x\pm by),ax\mp y\},\{-(x\pm by),-(ax\mp y)\}$.
 Since $x$ is odd, we may choose the sign of $x$
so that $x\e 1\mod 4$. Then $x\pm by\e (-1)^{\f{a-1}2+\f y2}(ax\mp
y)\e 1\mod 4$.
 Therefore, applying Lemma 2.5 we have
$$\aligned\la (a,b;n+1)&=\sum\Sb X\e Y\e 1\mod 4\\aX^2+bY^2=8n+a+b\endSb XY
\\&=(x+by)\cdot (-1)^{\f{a-1}2+\f y2}(ax-y)+(x-by)\cdot
(-1)^{\f{a-1}2+\f y2}(ax+y)\\&=(-1)^{\f{a-1}2+\f
y2}2(ax^2-by^2)=(-1)^{\f{a-1}2+\f y2}(4ax^2-2p)\\&=(-1)^{\f{a-1}2+\f
y2}(2p-4by^2).\endaligned$$  This together with Lemma 2.9 proves the
theorem.

\par\q\par
As examples, taking $a=3,5$ and $b=2$ in Theorem 3.3 we have:
$$\align&(-1)^{\f y2}(8y^2-2p)=\la(2,3;(7p+3)/8)\qtq{for}
p=3x^2+2y^2\e 11\mod{24},\tag 3.5
\\&(-1)^{\f y2}(2p-8y^2)=\la(2,5;(11p+1)/8)\qtq{for}p=5x^2+2y^2\e
13,37\mod{40}.\tag 3.6\endalign$$

\pro{Corollary 3.5} Let $a,b\in\Bbb N$ with $2\nmid a$, $2\mid b$,
$8\nmid b$ and $(a,b)=1$. Let $p \e 1\mod 8$ be a prime such that
$p\not=ab,ab+1$ and $p=x^2+aby^2$ with $x,y\in\Bbb Z$. Then
$$\lambda\Big(a,b;\f{(a+b)(p-1)}8+1\Big)=\lambda\Big(1,ab;\f{(ab+1)(p-1)}8+1\Big).$$
\endpro
Proof. By Theorem 3.3 we have
$$(-1)^{\f y2}(4x^2-2p)=\lambda\Big(1,ab;\f{(ab+1)(p-1)}8+1\Big).$$
This together with Theorem 3.2(ii) gives the result.

\subheading{4. Constructing $xy$ for primes $p=ax^2+by^2$}

\pro{Theorem 4.1} Let $a,b\in\Bbb N$, $8\nmid a$, $8\nmid b$ and
$n\in\{0,1,2,\ldots\}$. Let $p$ be an odd prime such that
$p=8n+a+b=ax^2+by^2$ with $x,y\in\Bbb Z$ and $x\e y\mod 4$. Then
$$xy=\la(a,b;n+1)\qtq{and}
2ax^2-p=\pm\sqrt{p^2-4ab\la(a,b;n+1)^2}.$$
\endpro
Proof. Let $x,y\in\Bbb Z$ be such that $p=8n+a+b=ax^2+by^2$. We
claim that $2\nmid xy$.  When $2\mid x$, we have $2\nmid y$, $a\e
8n+a= ax^2+by^2-b\e ax^2\e 0,4a\mod 8$ and so $8\mid a$. When $2\mid
y$, we have $2\nmid x$, $b \e 8n+b=ax^2+by^2-a\e by^2\e 0,4b\mod 8$
and so $8\mid b$. As $8\nmid a$ and $8\nmid b$, we see that $2\nmid
xy$. Suppose $x\e y\e 1\mod 4$. Then $x$ and $y$ are unique by Lemma
2.1. Now applying Lemma 2.5 we obtain $xy=\la(a,b;n+1)$.
\par Set $\la=\la(a,b;n+1)$. Then $x^2(p-ax^2)=bx^2y^2=b\la^2$ and
so $ax^4-px^2+b\la^2=0$. Thus, $x^2=(p\pm\sqrt{p^2-4ab\la^2})/(2a)$.
This completes the proof.

 \pro{Theorem 4.2} Let $a,b\in\Bbb N$ with $(a,b)=1$ and $ab\e 1\mod 4$.
 Let $p$ be an odd prime and $2p=8n+a+b=ax^2+by^2$ with
$n\in\{0,1,2,\ldots\}$, $x,y\in\Bbb Z$ and $4\mid x-y$. Then
$$xy=\la(a,b;n+1)\qtq{and}ax^2=p\pm\sqrt{p^2-ab\la(a,b;n+1)^2}.$$
\endpro
Proof. Let $x,y\in\Bbb Z$ be such that $2p=8n+a+b=ax^2+by^2$. Then
clearly $2\nmid xy$. Suppose $x\e y\e 1\mod 4$. Then $x$ and $y$ are
unique by Lemma 2.6. Now applying Lemma 2.5 we obtain $xy=\la$,
where $\la=\la(a,b;n+1)$. Thus, $x^2(2p-ax^2)=bx^2y^2=b\la^2$ and so
$ax^4-2px^2+b\la^2=0$. Hence, $x^2=(p\pm\sqrt{p^2-ab\la^2})/a$. This
completes the proof.

\pro{Theorem 4.3} Let $a,b\in\Bbb N$, $2\nmid ab$, $ab\not=3$,
$a+b\e 4\mod 8$. Let $p$ be an odd prime such that $p\nmid ab$ and
$4p=ax^2+by^2$ with $x,y\in\Bbb Z$ and $x\e y\e 1\mod 4$. Then
$$xy=\la\qtq{and}ax^2=2p\pm\sqrt{4p^2-ab\la^2},$$
where $\la=\la(a,b;\f 12(p-\f{a+b}4)+1)$.
\endpro
 Proof. Clearly $ab\e 3\mod
4$. From Lemma 2.7 we know that  $x$ and $y$ are unique. Set $n=\f
12(p-\f{a+b}4)$. Then $8n+a+b=4p$. By Lemma 2.5 we have $xy=\la$ and
so $x^2y^2=\la^2$. Thus $x^2(4p-ax^2)=b\la^2$ and so
$ax^4-4px^2+b\la^2=0$. Hence
$x^2=\f{4p\pm\sqrt{16p^2-4ab\la^2}}{2a}=\f{2p\pm
\sqrt{4p^2-ab\la^2}}a$. This completes the proof.

\par For example, if $p\not=11$ is an odd prime and
$4p=x^2+11y^2$ with $x\e y\e 1\mod 4$, then $xy=\la$ and
$x^2=2p\pm\sqrt{4p^2-11\la^2}$, where $\la=\la(1,11;(p-1)/2)$.

\subheading{5. Evaluation of $\la(1,3;n),\la(1,7;2n+1)$ and
$\la(3,5;2n+1)$}

\par\q For $n\in\Bbb N$, in [6, Vol.2, p.377] Klein and Fricke showed that
$$\la(1,1;n+1)=\sum\Sb x,y\in\Bbb Z,x\e 1\mod 4
\\x^2+y^2=4n+1\endSb (x^2-y^2).\tag 5.1$$
See also [8]. In the section we evaluate $\la(1,3;n),\la(1,7;2n+1)$
and $\la(3,5;2n+1)$.

\pro{Lemma 5.1} Let $a,b,n\in\Bbb N$ with $ab\e 3\mod 4$. Then
 $$\sum\Sb x+ay\e 1\mod 4\\x^2+aby^2=2n+1\endSb (x+ay)(x-by)=\f
 12\sum\Sb x,y\in\Bbb Z\\x^2+aby^2=2n+1\endSb (x^2-aby^2).$$
 \endpro
Proof. If $2\mid n$, then $x^2+aby^2=2n+1$ implies $x-1\e y\e 0\mod
2$. Thus
$$\aligned &\sum\Sb x,y\in\Bbb
Z,x+ay\e 1\mod 4
\\x^2+aby^2=2n+1\endSb (x+ay)(x-by)
\\&=\sum\Sb x\in\Bbb Z, x\e 1\mod 4\\x^2=2n+1\endSb x^2
+\sum\Sb x\e 1\mod 4,4\mid y, y>0\\x^2+aby^2=2n+1\endSb
\big\{(x+ay)(x-by)+(x-ay)(x+by)\big\} \\&\qq+ \sum\Sb x\e 3\mod
4,4\mid y-2, y>0\\x^2+aby^2=2n+1\endSb
\big\{(x+ay)(x-by)+(x-ay)(x+by)\big\}
\\&=\sum\Sb x\e 1\mod 4, y\in\Bbb Z\\x^2+aby^2=2n+1\endSb (x^2-aby^2)
=\f 12 \sum\Sb x, y\in\Bbb Z\\x^2+aby^2=2n+1\endSb (x^2-aby^2).
\endaligned$$ If $2\nmid n$, then $x^2+aby^2=2n+1$ implies
$x\e y+1\e 0\mod 2$. Thus
$$\aligned&\sum\Sb x,y\in\Bbb
Z,x+ay\e 1\mod 4
\\x^2+aby^2=2n+1\endSb (x+ay)(x-by)
\\&=\sum\Sb y\in\Bbb Z, y\e a\mod 4\\aby^2=2n+1\endSb (-aby^2)
+\sum\Sb 4\mid y+a,4\mid x-2, x>0\\x^2+aby^2=2n+1\endSb
\big\{(x+ay)(x-by)+(-x+ay)(-x-by)\big\} \\&\qq+ \sum\Sb y\e a\mod
4,4\mid x, x>0\\x^2+aby^2=2n+1\endSb
\big\{(x+ay)(x-by)+(-x+ay)(-x-by)\big\}
\\&=\sum\Sb y\e a\mod 4, x\in\Bbb Z\\x^2+aby^2=2n+1\endSb (x^2-aby^2)
=\f 12 \sum\Sb x, y\in\Bbb Z\\x^2+aby^2=2n+1\endSb (x^2-aby^2).
\endaligned$$
Thus the lemma is proved.

 \pro{Theorem 5.1} Let $n\in\Bbb N$. Then
$$\aligned&\la(1,3;n+1)=\f 12\sum\Sb x,y\in\Bbb Z
\\x^2+3y^2=2n+1\endSb (x^2-3y^2),
\\&\la(1,7;2n+1)=\f 12\sum\Sb x,y\in\Bbb Z
\\x^2+7y^2=2n+1\endSb (x^2-7y^2).\endaligned$$
\endpro
Proof. From Lemma 2.5 we have
$$\la(1,3;n+1)=\sum\Sb X\e Y\e 1\mod 4\\X^2+3Y^2=8n+4\endSb XY.$$
As $H(-12)=\{[1,0,3]\}$, by Lemma 2.2 we have
$$R([1,0,3],8n+4)=\f
12R([1,0,3],4)R([1,0,3],2n+1)=3R([1,0,3],2n+1).$$ Thus, if
$R([1,0,3],2n+1)=0$, then $R([1,0,3],8n+4)=0$ and so
$\la(1,1;n+1)=0$. Hence the result is true in this case.  Now assume
that $2n+1=x^2+3y^2$ with $x,y\in\Bbb Z$.
 Then $8n+4=4(x^2+3y^2)=(x+3y)^2+3(x-y)^2$.
As $R([1,0,3],8n+4)=3R([1,0,3],2n+1)$ we see that all the integral
solutions to the equation $8n+4=X^2+3Y^2$ are given by
$\{2x,2y\},\{x+3y,x-y\},\{x+3y,-(x-y)\}$, where $\{x,y\}$ runs over
all integral solutions to the equation $2n+1=x^2+3y^2$. Hence, using
Lemmas 2.5, 5.1 and the fact $x+3y\e x-y\mod 4$ we deduce
$$\aligned\la(1,3;n+1)&=\sum\Sb X\e Y\e 1\mod 4\\X^2+3Y^2=8n+4\endSb XY=\sum\Sb x,y\in\Bbb
Z,x+3y\e 1\mod 4
\\x^2+3y^2=2n+1\endSb (x+3y)(x-y)
\\&=\f 12\sum\Sb x,y\in\Bbb Z\\x^2+3y^2=2n+1\endSb (x^2-3y^2).\endaligned$$

\par Now we consider the formula for $\la(1,7;2n+1)$. By Lemma 2.5 we have
$$\la(1,7;2n+1)=\sum\Sb X\e Y\e 1\mod 4\\X^2+7Y^2=16n+8\endSb XY.$$
As $H(-28)=\{[1,0,7]\}$, by Lemma 2.2 we have
$$R([1,0,7],16n+8)=\f
12R([1,0,7],8)R([1,0,7],2n+1)=2R([1,0,7],2n+1).$$ Thus, if
$R([1,0,7],2n+1)=0$, then $R([1,0,7],8(2n+1))=0$ and so
$\la(1,7;2n+1)=0$. Hence the result is true in this case.  Now
assume that $2n+1=x^2+7y^2$ with $x,y\in\Bbb Z$.
 Then $16n+8=8(x^2+7y^2)=(x+7y)^2+7(x-y)^2$.
As $R([1,0,7],16n+8)=2R([1,0,7],2n+1)$ we see that all the integral
solutions to the equation $16n+8=X^2+7Y^2$ are given by
$\{x+7y,x-y\},\{x+7y,-(x-y)\}$, where $\{x,y\}$ runs over all
integral solutions to the equation $2n+1=x^2+7y^2$. Hence, using
Lemmas 2.5, 5.1 and the fact $x+7y\e x-y\mod 4$ we deduce
$$\aligned\la(1,7;2n+1)&=\sum\Sb X\e Y\e 1\mod 4\\X^2+7Y^2=16n+8\endSb XY
=\sum\Sb x,y\in\Bbb Z,x+7y\e 1\mod 4
\\x^2+7y^2=2n+1\endSb (x+7y)(x-y)
\\&=\f 12\sum\Sb x,y\in\Bbb Z\\x^2+7y^2=2n+1\endSb (x^2-7y^2).\endaligned$$
This completes the proof.
 \pro{Theorem 5.2} For $n\in\Bbb N$ we have
$$\la(1,15;4n+1)=\la(3,5;2n+1)=\f 12\sum\Sb x,y\in\Bbb
Z\\x^2+15y^2=2 n+1\endSb (x^2-15y^2).$$
\endpro
Proof.  It is clear that
$$R([3,0,5],8)=4,\ R([1,0,15],8)=0,\ R([1,0,15],16)=6\qtq{and}
R([3,0,5],16)=0.$$ As $H(-60)=\{[1,0,15],[3,0,5]\}$, by Lemma 2.2
and the above we have
$$\align &2R([3,0,5],8(2n+1))\\&=R([3,0,5],8)R([1,0,15],2n+1)+R([1,0,15],8)R([3,0,5],2n+1)
\\&=4R([1,0,15],2n+1)\endalign$$
and
$$\align &2R([1,0,15],16(2n+1))\\&=R([1,0,15],16)R([1,0,15],2n+1)
+R([3,0,5],16)R([3,0,5],2n+1)
\\&=6R([1,0,15],2n+1).\endalign$$
From Lemma 2.5 we have
$$\la(1,15;4n+1)=\sum\Sb X\e Y\e 1\mod 4\\X^2+15Y^2=16(2n+1)\endSb
XY\qtq{and}\la(3,5;2n+1)=\sum\Sb X\e Y\e 1\mod
4\\3X^2+5Y^2=8(2n+1)\endSb XY.$$ Thus, if $R([1,0,15],2n+1)=0$, then
$R([1,0,15],16(2n+1))=R([3,0,5],8(2n+1))=0$ and so
$\la(1,15;4n+1)=\la(3,5;2n+1)=0$. Hence the result is true in this
case.
\par Now assume that $2n+1=x^2+15y^2$ with $x,y\in\Bbb Z$. Then
$16(2n+1)=(4x)^2+15(4y)^2=(x+15y)^2+15(x-y)^2$ and
$8(2n+1)=3(x+5y)^2+5(x-3y)^2$. Since $R([3,0,5],
8(2n+1))=2R([1,0,15],2n+1)$ and
$R([1,0,15],16(2n+1))=3R([1,0,15],2n+1)$, we see that all the
integral solutions to the equation $3X^2+5Y^2=8(2n+1)$ are given by
$\{x+5y,\pm (x-3y)\}$, and  all the integral solutions to the
equation $X^2+15Y^2=16(2n+1)$ are given by $\{4x,4y\}$ and
$\{x+15y,\pm (x-y)\}$, where $\{x,y\}$ runs over all integral
solutions to the equation $2n+1=x^2+15y^2$. As $x+5y\e x-3y\e \pm
1\mod 4$ and $x+15y\e x-y\e \pm 1\mod 4$, from the above and Lemma
5.1 we deduce
$$\aligned \la(1,15;4n+1)&=\sum\Sb X\e Y\e 1\mod 4\\X^2+15Y^2=16(2n+1)\endSb
XY=\sum\Sb x+15y\e 1\mod 4\\x^2+15y^2=2n+1\endSb (x+15y)(x-y)
\\&=\f 12\sum\Sb x,y\in\Bbb Z\\x^2+15y^2=2n+1\endSb (x^2-15y^2)\endaligned$$
and
$$\aligned \la(3,5;2n+1)&=\sum\Sb X\e Y\e 1\mod 4\\3X^2+5Y^2=8(2n+1)\endSb
XY=\sum\Sb x+5y\e 1\mod 4\\x^2+15y^2=2n+1\endSb (x+5y)(x-3y)
\\&=\f 12\sum\Sb x,y\in\Bbb Z\\x^2+15y^2=2n+1\endSb
(x^2-15y^2).\endaligned$$ This proves the theorem.

 \pro{Theorem 5.3} Let $p>5$ be a prime.
Then
$$\aligned &\la(3,5;p)=\cases 0&\t{if $p\not\e 1,19\mod{30}$,}
\\4x^2-2p&\t{if $p\e 1,19\mod{30}$ and so $p=x^2+15y^2(x,y\in\Bbb Z)$,}\endcases
\\&\la(3,5;2p)=\cases 0&\t{if $p\not\e 17,23\mod{30}$,}
\\2p-12x^2&\t{if $p\e 17,23\mod{30}$ and so $p=3x^2+5y^2(x,y\in\Bbb Z)$,}\endcases
\\&\la(3,5;3p)=\cases 0&\t{if $p\not\e 17,23\mod{30}$,}
\\36x^2-6p&\t{if $p\e 17,23\mod{30}$ and so $p=3x^2+5y^2(x,y\in\Bbb Z)$,}\endcases
\\&\la(3,5;5p)=\cases 0&\t{if $p\not\e 17,23\mod{30}$,}
\\10p-60x^2&\t{if $p\e 17,23\mod{30}$ and so $p=3x^2+5y^2(x,y\in\Bbb Z)$.}
\endcases
\endaligned$$
\endpro
Proof. If $p\e 1,19\mod{30}$, then $p=x^2+15y^2$ for some positive
integers $x$ and $y$ (see [14, Table 9.1]). By Lemma 2.1, $x$ and
$y$ are unique. From Theorem 5.2 we have
$$\la(3,5;p)=\f 12\sum\Sb x,y\in\Bbb
Z\\x^2+15y^2=p\endSb (x^2-15y^2)=2(x^2-15y^2)=4x^2-2p.$$
 If $p\not\e 1,19\mod{30}$, then $p$ is not
represented by $x^2+15y^2$. Thus, by Theorem 5.2 we have
$\la(3,5;p)=0$.

\par If $p\e 17,23\mod{30}$, then $p=3x^2+5y^2$ with
$x,y\in\Bbb Z$. Taking $a=3$ and $b=5$ in Theorem 3.1 we obtain
$\la(3,5;2p)=2p-12x^2$. If $p\not\e 17,23\mod{30}$, as
$R([3,0,5],2)=R([1,0,15],2)=0$, using Lemma 2.2 we see that
$2R([3,0,5],2p)=R([3,0,5],2)R([1,0,15],p)+R([1,0,15],2)R([3,0,5],p)=0$.
Thus, appealing to Lemma 2.5 we have $\la(3,5;2p)=0$.
\par Let $b\in\{3,5\}$. By Theorem 5.2 we have
$$\la(3,5;bp)=\f 12\sum\Sb X,Y\in\Bbb Z\\X^2+15Y^2=bp\endSb
(X^2-15Y^2).$$ As $H(-60)=\{[1,0,15],[3,0,5]\}$, $R([1,0,15],b)=0$
and $R([3,0,5],b)=2$, using Lemma 2.2 we see that
$$\align R([1,0,15],bp)&=\f
12(R([1,0,15],3)R([1,0,15],p)+R([3,0,5],3)R([3,0,5],p))
\\&=R([3,0,5],p).\endalign$$
If $p\not\e 17,23\mod{30}$, then
$R([1,0,15],bp)=R([3,0,5],p)=0$ and so $\la(3,5;bp)=0$. If $p\e
17,23\mod{30}$, then there are unique positive integers $x$ and $y$
such that $p=3x^2+5y^2$. As $R([1,0,15],bp)=R([3,0,5],p)=4$, we see
that all the integral solutions to $3p=X^2+15Y^2$ are given by
$\{\pm 3x,\pm y\}$, and all the integral solutions to $5p=X^2+15Y^2$
are given by $\{\pm 5y,\pm x\}$. Thus,
$$\la(3,5;3p)=\f 12\sum\Sb X,Y\in\Bbb Z\\X^2+15Y^2=3p\endSb
(X^2-15Y^2)=2((3x)^2-15y^2)=36x^2-6p$$ and
$$\la(3,5;5p)=\f 12\sum\Sb X,Y\in\Bbb Z\\X^2+15Y^2=5p\endSb
(X^2-15Y^2)=2((5y)^2-15x^2)=10p-60x^2.$$ This completes the proof.

\enddocument
\bye